\documentclass[10pt,a4paper]{amsart}
\usepackage{amsfonts}
\usepackage{amsmath}
\usepackage{amscd,amssymb,amsthm}
\usepackage{graphics}

\begin{document}
\title{Analysis of the reconnection process in nontwist cubic maps }
\author{Gheorghe Tigan}
\maketitle

\begin{abstract}
The reconnection process in the dynamics of cubic nontwist maps,
introduced in [3], is studied. The present paper extends the work
presented in [8]. As in that work, in order to describe the route
to reconnection of the involved Poincar\'{e}--Birkhoff chains or
dimerised chains we investigate an approximate interpolating
Hamiltonian of the map under study revealing again that the
scenario of reconnection of cubic nontwist maps is different from
that occurring in the dynamics of quadratic nontwist maps.
\end{abstract}

\textit{Key words}: area preserving maps,  nontwist maps,
reconnection bifurcation.\newline \

\emph{email}: gtigan73@yahoo.com

\section{Introduction}

Nontwist maps arise naturally in the study of Hamiltonian systems,
because they are models for Poincar\'{e} maps associated to
sections in an energy manifold of an iso--energetically degenerate
two degree of freedom Hamiltonian system \cite{cast96, voyatzis,
apte}, in transport problems in plasma physics, accelerator
physics and in other areas. Transport problems in plasma physics
can be modelled by an area preserving map where the twist
condition fails \cite{doyle}. Applications of nontwist maps in
accelerator physics can be found in \cite{gerasimov}. In
\cite{petrisor21} are studied the quadratic nontwist standard-like
maps both from theoretical and numerical point of view. Nontwist
standard-like maps exhibit both time-reversal and spatial symmetry
being observed the appearance of the meanders. Meanders are
invariant circles that exhibit foldings in such a way that they
are not graphs of functions. During the last decade numerical and
theoretical studies of quadratic non--twist maps \cite{howard95},
\cite{cast96}, \cite{petrisor21}, \cite{petrisor22} revealed a
global bifurcation (called reconnection) of the invariant
manifolds of two distinct regular hyperbolic periodic orbits
having the same rotation number. At the threshold of reconnection
the involved hyperbolic orbits are connected by a common arc of
their invariant manifolds. The physical model of reconnection is
met in Tokamaks \cite{viana} which are experimental machines for
achievement of controlled thermonuclear fusion reactions.
\par For a rigorous analysis of local and global bifurcations
occurring in a family of area preserving maps defined on an
annulus $\mathbb{T}\times [a,b]$ ($\mathbb{T}$ denotes the circle
identified with $[0,2\pi)$) one derives an approximate
interpolating Hamiltonian of the map under study \cite{simo}.

\par
The present work deals with reconnection in  the cubic nontwist
area preserving diffeomorphism of the annulus
$\mathbb{T}\times\mathbb{R}$, $f:(x,y)\mapsto (x',y')$:
\begin{equation}\label{cubntw}\begin{array}{lll}
x'&=&x+F(a,b;y')\,\,(\mbox{mod}\,\, 2\pi)\\
y'&=&y+k\sin{x}
\end{array}\end{equation}
where the rotation number function $F$ is a cubic map depending on
two parameters $a>0, b\in\mathbb{R}, b\neq 0$, i.e.
$F(a,b;y)=y-ay^{2}+by^{3}$ and $k>0$ is a perturbation
parameter.\par We recall that an area preserving diffeomorphism
$g:\mathbb{T}\times\mathbb{R}\to\mathbb{T}\times\mathbb{R}$,
$g:(x,y)\mapsto(x',y')$ is a twist map if $\partial_y x'\neq 0$
($\partial_y$ denotes the partial derivative with respect to $y$).
 Twist property is a basic assumption of KAM theorem,
as well as for the Aubry--Mather theory \cite{katok}. The map (
\ref{cubntw}) is a non--twist map because it violates the twist
condition. Our purpose is to study its dynamics as well as the
route to reconnection in the case when the shape parameter $a$ and
the perturbation parameter $k$ are fixed and the other shape
parameter $b$ varies on the real line or on an interval.

\section{Properties of cubic nontwist map}

First we recall some properties of the map under study
\cite{tig1}. The motion in the unperturbed map (\ref{cubntw}), i.e
the map corresponding to $k=0$:
\begin{equation}\label{intntw}\begin{array}{lll}
x'&=&x+y-ay^{2}+by^{3}\,\,(\mbox{mod}\,\, 2\pi)\\
y'&=&y\end{array}\end{equation} occurs along the circle
$y=\mbox{cst}$. The rotation number of an orbit starting at
$(x,y)$ is:
\begin{equation}
\rho=\lim_{n\to\infty}\displaystyle\frac{X_n-X}{2n\pi}=F(a,b;y)/(2\pi),
\end{equation}
where $(X_n,Y_n)$ is the orbit of the point $(X,Y)=(x,y)$ under
the lift of the map (the map defined on $\mathbb{R}^2$ having the
same expression, without modulo $2\pi$ for the first component).
The map (\ref{intntw}) violates the twist condition for the
parameter values $(a,b)$ such that $a^2-3b\geq 0$, along the
circles:
\begin{equation}
C_1:\,\, y=\displaystyle\frac{a+\sqrt{a^2-3b}}{3b},\quad
C_2:\,\,y=\displaystyle\frac{a-\sqrt{a^2-3b}}{3b}\end{equation}

\noindent These circles are called {\it twistless} or {\it
shearrless} circles. At the same time along the circle $C_1$ the
rotation number has a global minimum,

$F|_{C_1} =\frac{-2a^{3}-3a^{2}\sqrt{a^{2}-3b}+9ab+(a^{2}+6b)\sqrt{a^{2}-3b}}{%
54\pi b^{2}}$,

while along $C_2$, a global maximum,

$F|_{C_2} =\frac{-2a^{3}+3a^{2}\sqrt{a^{2}-3b}+9ab-(a^{2}+6b)\sqrt{a^{2}-3b}}{%
54\pi b^{2}}$.

\bigskip

\par
\begin{figure}[h]
\includegraphics{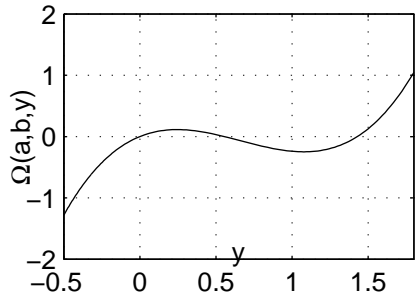}
\caption{The graph of the rotation number function $F(a,b;y)$ for
$a=2.5$ and $b=1.26$ } \label{fig1b}
\end{figure}

 Let us denote by $\omega_m, \omega_M$:
\begin{equation}
\omega_m:=\displaystyle\frac{a+\sqrt{a^2-3b}}{3b},\quad
\omega_M:=\displaystyle\frac{a-\sqrt{a^2-3b}}{3b},\end{equation}
the points of minimum, respectively maximum, for the rotation
number function $F$. For $y\in(-\infty, \omega_M )$ the
unperturbed map has positive twist (the rotation number function
is increasing), for $y\in(\omega_M,\omega_m)$ has a negative twist
(the rotation number function is decreasing), while for $y\in
(\omega_m,+\infty)$ has again a positive twist (Fig.\ref{fig1b}).

The orbits lying on the circles $y=y_0$ with $F(a,b;y_0)/2\pi$ a
rational number in lowest terms, $p/q$, are periodic orbits. If
such a periodic orbit lies in a region of monotone twist property
of the map, after a slight perturbation, it gives rise generically
to at least two periodic orbits of the same rotation number, one
elliptic and the second regular hyperbolic. Elliptic points are
surrounded by invariant circles, and hyperbolic points are
connected by heteroclinic connections. Such a pair of periodic
orbits and the associated invariant sets form a
Poincar\'{e}--Birkhoff chain. For $y\in (\omega_M,\omega_m)$ can
exist three circles $y=\mbox{cst}$ of the unperturbed map, on
which lie periodic orbits of the same rotation number $p/q$. Our
aim is to study the bifurcations of the periodic orbits or the
invariant manifolds belonging to three distinct
Poincar\'{e}--Birkhoff chains created after a slight perturbation,
as the shape parameter $b\in\mathbb{R}$, $a^2-3b>0$ defining the
rotation number function $F$ varies.

\section{Reconnection scenario}

In order to analyze the changes in the topology of invariant
manifolds of the involved p/q-type hyperbolic periodic orbits
consider the interpolating Hamiltonian associated to the map $F$:

\begin{equation}H_{a,b,k}(x,y)=-y^{2}/2+ay^{3}/3-by^{4}/4-k\cos x
\end{equation}
It defines the vector field
\begin{equation} X_{H}=\left( \frac{-\partial H}{\partial y},\frac{\partial H}{\partial x}%
\right) =\left( y-ay^{2}+by^{3},k\sin x\right) \end{equation}
which is reversible with respect to the involution $R(x,y)=(-x,y),i.e.$ \ $%
R\circ X_{H}=-X_{H}\circ R.$ The fixed point set, $\mbox{Fix}(R)$,
consists in the lines $x=0$ and $x=\pi ,$ called symmetry lines.
The equilibrium \ points of \ $X_{H}$ lying on the symmetry lines
are called symmetric. The Hamiltonian system associated to the
vector field $X_{H}$ can display at most three chains:
Poincar\'{e}--Birkhoff chains or dimerised chains. A dimerised
chain is a structure formed by elliptic points surrounded by
homoclinic circles to the corresponding hyperbolic points. \par In
order to describe the scenario of reconnection and the local
bifurcations of the equilibrium points, we analyze the position on
the symmetry lines of the equilibrium points, their stability type
and bifurcations occurring as $b$ varies and $a$, $k$ are fixed in
the parametric space $(a,b,k)$, with $a,k>0$ and $b\in\mathbb{R}$.
Therefore, if $b<\frac{a^{2}}{4}$, the Hamiltonian system has six
equilibrium points (e stands for elliptic and h for hyperbolic):
\bigskip

$P_{1h}(0,0);P_{2e}(\pi
,0);P_{3e}(0,\frac{a-\sqrt{a^{2}-4b}}{2b});$

\bigskip

$P_{4h}(\pi ,\frac{a-\sqrt{a^{2}-4b}}{2b});P_{5h}(0,\frac{a+\sqrt{a^{2}-4b}}{%
2b});P_{6e}(\pi ,\frac{a+\sqrt{a^{2}-4b}}{2b})$

\bigskip

\begin{figure}[h]
\includegraphics{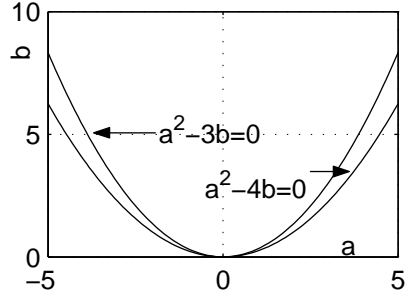}
\caption{The existence domain of the equilibrium points }
\label{parb}
\end{figure}

If $\frac{a^{2}}{4}<b<\frac{a^{2}}{3}$ (Fig.\ref{parb}) the vector
field $X_{H}$ has only two equilibrium points: $ P_{1}(0,0),
P_{2}(\pi ,0)$, while for $a^{2}=4b$ it has four equilibrium
points: $ P_{1}(0,0), P_{2}(\pi ,0), A(0,\frac{1}{\sqrt{b}}),
B(\pi,\frac{1}{\sqrt{b}})$. In the latter case, $a^2-4b=0$, the
two eigenvalues are zero and a bifurcation of equilibrium points
occurs.

\bigskip

\begin{figure}[h]
\includegraphics{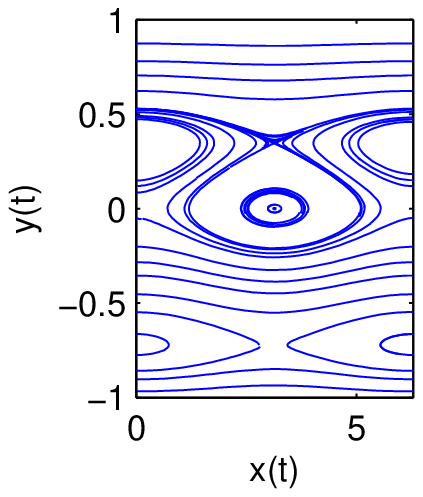}
\includegraphics{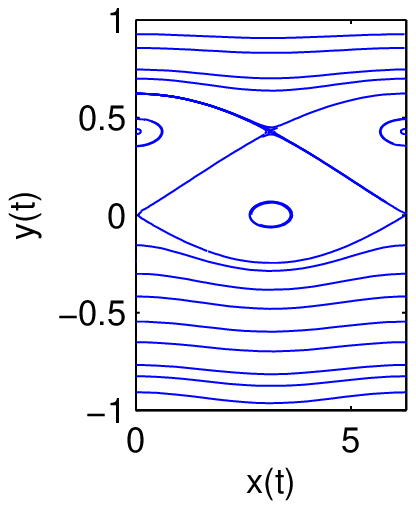}
\includegraphics{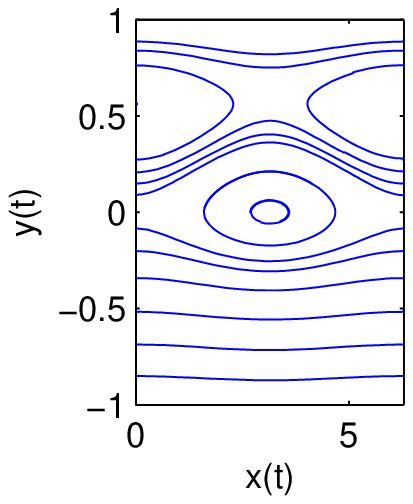}
\caption{ Reconnection scenario of the chains I and II. The values
of the parameter $b$ are: a)\ $b=-4$ (left) \ b)\ $b=-1.9538$
(middle) \ c)\  $b=-0.5$ (right) } \label{recbjos}
\end{figure}

In the following we want to describe the local changes in the
topology of the invariant manifolds of the Hamiltonian system,
when $a,k> 0 $ are fixed and $b\neq 0$ varies on the real line. We
remark that to get connected any two neighboring chains when $b$
varies, we need to consider both positive and negative values of
the parameter $b$. When the parameter $a$ varies \cite{tig1}, it
is sufficient to consider the case $a>0$. The systems whose phase
portraits are illustrated in different figures correspond to
$a=1.5$ and $k = 0.018$. Denote by I, II and III the three chains
containing the equilibrium points, more precisely, the chain I
contains the points $P_{1h}, P_{2e}$, the chain II contains the
points $P_{3e}, P_{4h}$ and III the points $P_{5h}, P_{6e}$. For
$b$ small enough ($b<-2$ for example), the all six equilibrium
points are born, Fig.\ref{recbjos}a). The chains I and II are two
dimerised chains. Between these chains the trajectories of the
Hamiltonian vector field are not graphs of real functions of $x$,
but they are \emph{meanders}. Each two neighboring points on the
same symmetry line have opposite stability type. Increasing
further the parameter $b$, the equilibrium points lying on the
same symmetry line (such points lie within two different chains)
go away. At a critical value, called threshold of reconnection,
the hyperbolic points of the two neighboring chains get connected
by common branches of their invariant manifolds,
Fig.\ref{recbjos}b). The common branches are $W^{u}(P_{1h})$ and
$W^{s}(P_{4h})$ respectively $W^{s}(P_{1h})$ and $W^{u}(P_{4h})$.
In order to get this threshold of reconnection we impose that the
hyperbolic equilibrium points $P_{1h}$ and $P_{4h}$ to belong to
the same energy level set, that is, $H_{a,b,k}(P_{1h})
=H_{a,b,k}(P_{4h})$. This implies that the reconnection surface of
the dimerised I and II chains is:
\begin{equation}\label{recc1}
6b^{2}+a^{4}-6a^{2}b+48b^{3}k+4ab\sqrt{a^{2}-4b}-a^{3}\sqrt{a^{2}-4b}=0
\end{equation}
Numerically, it leads to the first threshold of reconnection
$b:=b_{1rec}=-1.9538$.
\par
\begin{figure}[h]
\includegraphics{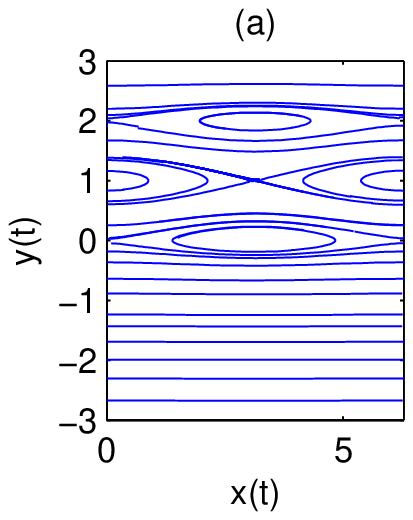}
\includegraphics{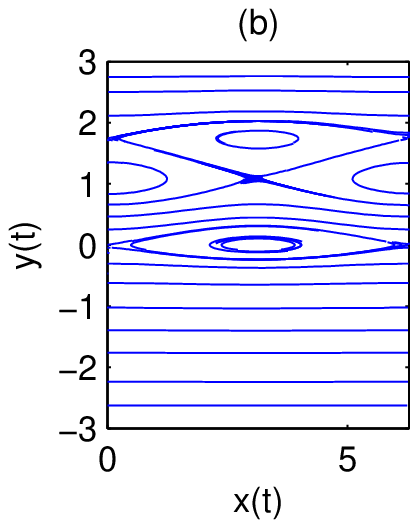}
\includegraphics{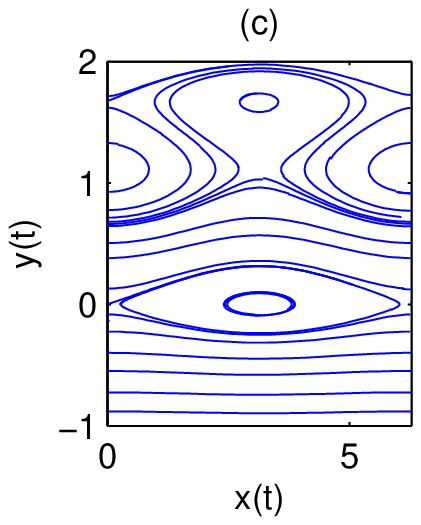}
\includegraphics{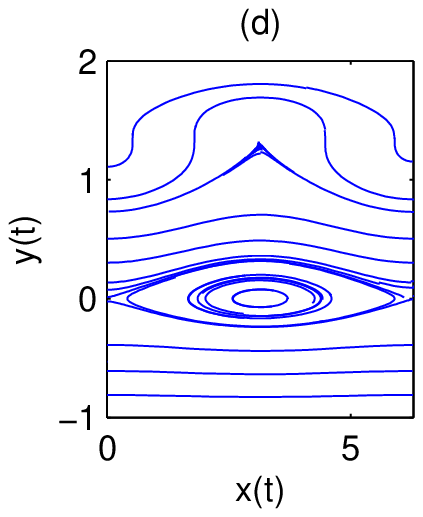}
\includegraphics{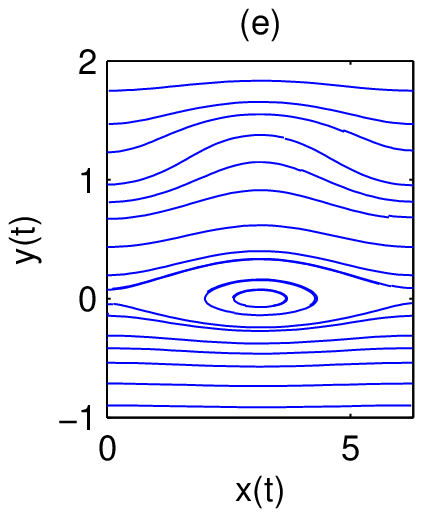}
\caption{ Reconnection scenario of the chains II and III. The
values of the parameter $b$ are: a)  $b=0.5$ \  b) $b=0.53168$ \
 c)  $b=0.54$ \  d)  $b=0.5625$ \  e)  $b=0.6$ } \label{recbsus}
\end{figure}

 Increasing $b$ slightly from $b_{1rec}$, the two dimerised
chains become two Poincar\'{e}--Birkhoff chains,
Fig.\ref{recbjos}c), so the system displays three
Poincar\'{e}--Birkhoff distinct chains, Fig.\ref{recbsus}a) and
the Poincar\'{e}--Birkhoff chains II and III approach each other.
At the threshold of reconnection the hyperbolic points of these
two chains (II and III) get connected by common branches of
invariant manifolds Fig.\ref{recbsus}b). The common branches are
$W^{u}(P_{4h})$ and $W^{s}(P_{5h})$ respectively $W^{s}(P_{4h})$
and $W^{u}(P_{5h})$. As above, imposing $H_{a,b,k}(P_{4h})
=H_{a,b,k}(P_{5h})$ we find the surface of reconnection of the
chains II and III:
\begin{equation}\label{recc2}
k=\frac{1}{24}\frac{a}{b^{3}}({a^{2}-4b})\sqrt{a^{2}-4b}
\end{equation}
For the numerical values, the threshold of reconnection is
$b:=b_{2rec}=0.53168$, Fig.\ref{recbsus}b). Continuing to increase
$b$ beyond $b_{2rec}$, the Poincar\'{e}--Birkhoff chains II and
III are transformed into two dimerised chains,
Fig.\ref{recbsus}c). Between these chains, the same as above for
the dimerised I and II chains, the trajectories of the Hamiltonian
vector field are meanders. At the value $b=\frac{a^{2}}{4}$, the
equilibrium points $P_3, P_4, P_5$ and $P_6$ are reduced to two
points $A(0,\frac{1}{\sqrt{b}}), B(\pi,\frac{1}{\sqrt{b}})$, \
Fig.\ref{recbsus}d),  which completely disappear for
$\frac{a^{2}}{4}<b<\frac{a^{2}}{3}$, (remember that we work only
on the nontwist domain, $b<\frac{a^{2}}{3}$), Fig.\ref{recbsus}e).

\bigskip

\bigskip

\bigskip

\begin{figure}[h]
\includegraphics{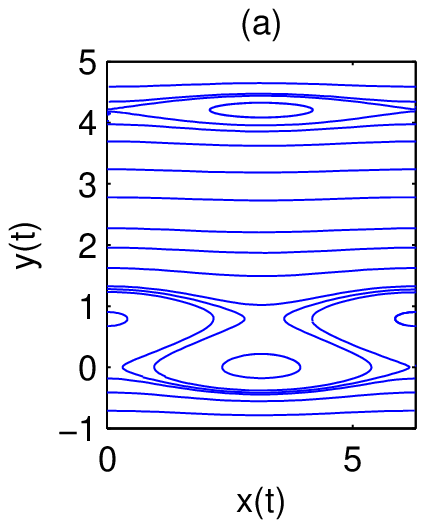}
\includegraphics{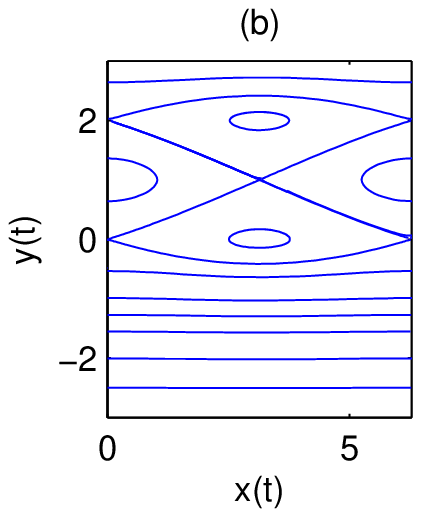}
\includegraphics{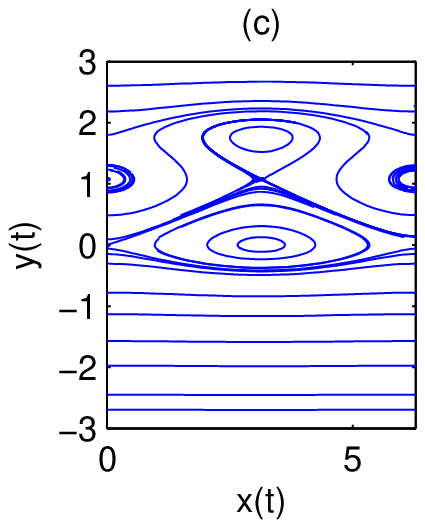}
\includegraphics{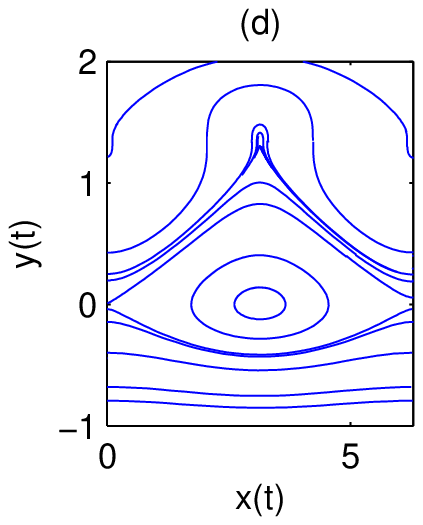}
\includegraphics{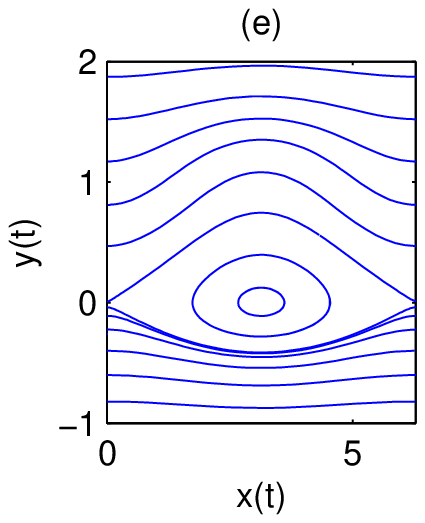}
\caption{ Reconnection scenario of the three chains I, II and III.
The values of the parameter $b$ are: a) $b=0.3$ \  b) $b=0.5$ \ c)
$b=0.53$ \  d)  $b=0.5625$ \  e)  $b=0.7$ \ } \label{doubleb}
\end{figure}
\bigskip
Remark that we can simultaneously get connected the all three
chains I, II and III. Call this \emph{triple reconnection}.  Let
us describe the triple reconnection scenario. Consider in this
case $b>0$. For $b$ slightly beyond $0$, the chains I and II are
two dimerised chains while III is a Poincar\'{e}--Birkhoff chain,
Fig.\ref{doubleb}a). Increasing further the parameter $b$, the
equilibrium points of the two dimerised chains lying on the same
symmetry line go away  while the equilibrium points of the
Poincar\'{e}--Birkhoff chain approaches the points of the
dimerised chain II. At a critical value, called triple threshold
of reconnection, the hyperbolic points of the three chains get
connected by common branches of their invariant manifolds
Fig.\ref{recbjos}b). Imposing $H_{a,b,k}(P_{1h}) =
H_{a,b,k}(P_{4h})=H_{a,b,k}(P_{5h}) $, we get from (\ref{recc1})
and (\ref{recc2}) the following reconnection curve:

\begin{equation}\label{triplerecb}
a^4 -6a^{2}b+6b^2+a(a^{2}-4b)\sqrt{a^{2}-4b}=0 \end{equation}

Solving numerically (\ref{triplerecb}) for $a=1.5$ one get
$b:=b_{3rec}=0.5$ and from (\ref{recc1}) $k:=k_{3rec}=0.0625$.
Consequently, to get connected the three chains I, II and III, in
the case when $b$ varies, we have to keep the perturbation
parameter $k$ at the constant value $k_{3rec}$. For $b>b_{3rec}$
the scenario is similar to the case $b>b_{2rec}$,
Fig.\ref{doubleb}c)-e) and Fig.\ref{recbsus}c)-e)

\section{Conclusions}

In this paper we have extended the studies reported in \cite{tig1}
on reconnection scenario of a three-parameter cubic nontwist map
depending on the parameters $a,b$ and $k$. Using an approximate
interpolating Hamiltonian of the map we have described the
reconnection process of any two neighboring chains in the case
when the parameters $a,k$ are fixed and $b$ varies. By numerically
computations we found the exact values of the thresholds of
reconnection. At the end we presented the triple reconnection of
the all three involved chains.

\section{Acknowledgements}

This work was (partially) supported through a European Community
Marie Curie Fellowship and in the framework of the CTS, contract
number HPMT-CT-2001-00278.


\bigskip



\bigskip

\begin{thebibliography}{50}

\bibitem{katok} A. Katok, B.  Hasselblatt, {\it Introduction to the modern theory
of dynamical systems}, Cambridge University Press, 1995.

\bibitem{cast96} D. del-Castillo-Negrete, J.M. Greene,  P.J. Morrison, {\it Area
preserving nontwist maps: periodic orbits and transition to
chaos}, Physica  D91 (1996) 1-23.

\bibitem{howard95} J. E. Howard, J. Humpherys, {\it Nonmonotonic twist maps},
Physica D80 (1995) 256--276.

\bibitem{petrisor21} E. Petrisor,  {\it Reconnection scenarios and the threshold
of reconnection in the dynamics of nontwist maps}, Chaos, Solitons
and Fractals, 14 (2002) 117-127.

\bibitem{petrisor22} E. Petrisor,  {\it Nontwist area preserving maps with
reversing symmetry group}, Int. J. Bif. Chaos, 11 (2001) 497-511.

\bibitem{simo} C. Sim\'{o},
 {\it Invariant curves of analytic perturbed nontwist area preserving
 maps}, Regular and Chaotic Dynamics, 3  (1998) 180-195.

\bibitem{doyle} E.J. Doyle et al., {\it Modifications in turbulence and edge
electric fields at the L–H transition in the DIII-D tokamak},
Physics of Fluids B Vol.3(8), (1991), 2300-2307.

\bibitem{tig1} Gh. Tigan, {\it On the scenario of reconnection in nontwist cubic maps,
Chaos, Soliton and Fractals (accepted), to appear.}

\bibitem{gerasimov} A. Gerasimov, F.M. Israilev, J.L. Tennyson,
A.B. Temnykh, {\it Springer Lectures Notes in Physics}, Vol.247,
(154), 1986.

\bibitem{soskin} S.M Soskin, {\it Phys.Rev.}, E 50(1), (1994), R44.

\bibitem{egydio} R. Egydio de Carvalho and A.M. Ozorio de Almeida,
{\it Integrable approximation to the overlap of resonances}, Phys.
Letters A, (1992), 162, 457-63.

\bibitem{chapel}  H.W. Chapel and T. Post, {\it The birth process of periodic
orbits in nontwist maps}, Physica D, (1995), 80, 256-276.

\bibitem{viana} R.L. Viana, {\it Chaotic magnetic field lines in a Tokamak with
resonant helical windings}, Chaos 11, (2000), 765-778.

\bibitem{voyatzis} G. Voyatzis, S. Ichtiaroglou, {\it Degenerate bifurcations of resonant tori in
     Hamiltonian systems}, Int. J. Bif. Chaos 9, (1999), 849--863.

 \bibitem{apte} A. Apte, A. Wurm, P.J. Morrison,  {\it Renormalization and destruction
of $1/\gamma^2$ tori in the standard nontwist map,} Chaos 13
(2003).

\end{thebibliography}
\end{document}